\documentclass[12pt,reqno]{amsart}
\usepackage{enumerate}
\usepackage{amsmath, amsfonts, amssymb, amsthm}
\def\sgn{{\rm sgn}\,}
\def\per{{\rm per}\,}

\numberwithin{equation}{section}
\newtheorem{Theorem} {Theorem} [section]
\newtheorem{Proposition} [Theorem] {Proposition}
\newtheorem{Lemma} [Theorem] {Lemma}
\newtheorem{Corollary} [Theorem] {Corollary}
\newtheorem{Conjecture}[Theorem]{Conjecture}
\theoremstyle{definition}

\newtheorem{Definition}[Theorem]{Definition}
\newtheorem{Remark}[Theorem]{Remark}

\begin{document}
\hbox{Israel J. Math. 182(2011), 425--437.}
\medskip
\title
[Exterior Algebras and Two Conjectures on Finite Abelian
Groups]{Exterior Algebras and Two Conjectures on Finite Abelian
Groups}
\author[Tao Feng, Zhi-Wei Sun and Qing Xiang]
{Tao Feng$^{1}$, Zhi-Wei Sun$^{2}$ and Qing Xiang$^{3}$}

\thanks{$^{1}$Research supported in part by the National Natural Science Foundation of China (grant 10331030).}
\thanks{$^{2}$Corresponding author. Supported by the National Natural Science Foundation (grant 10871087)
and the Overseas Cooperation Fund (grant 10928101) of China.}
\thanks{$^{3}$Research supported by the NSF Grant DMS 0701049 of USA and the Overseas Cooperation Fund (grant 10928101) of China.}
\address{Department of Mathematical Sciences, University of Delaware, Newark, DE 19716, USA}
\email{feng@math.udel.edu}
\address{Department of Mathematics, Nanjing University, Nanjing 210093, People's Republic of China} \email{zwsun@nju.edu.cn}
\address{Department of Mathematical Sciences, University of Delaware, Newark, DE 19716, USA} \email{xiang@math.udel.edu}


\begin{abstract} Let $G$ be a finite abelian group with $|G|>1$.
Let $a_1,\ldots,a_k$ be $k$ distinct elements of $G$ and let
$b_1,\ldots,b_k$ be (not necessarily distinct) elements of $G$,
where $k$ is a positive integer smaller than the least prime divisor
of $|G|$. We show that there is a permutation $\pi$ on
$\{1,\ldots,k\}$ such that $a_1b_{\pi(1)},\ldots,a_kb_{\pi(k)}$ are
distinct, provided that any other prime divisor of $|G|$ (if there
is any) is greater than $k!$. This in particular confirms the
Dasgupta-K\'arolyi-Serra-Szegedy conjecture for abelian $p$-groups.
We also pose a new conjecture involving determinants and characters,
and show that its validity implies Snevily's conjecture for abelian
groups of odd order. Our methods involve exterior algebras and
characters.
\end{abstract}

\maketitle

\section{Introduction}

Let $G=\{a_1,\ldots,a_n\}$ be an abelian group (written
multiplicatively) of order $n$, and let $b_1,\ldots,b_n\in G$. In
1952, M. Hall, Jr. \cite{hall} showed that
$a_1b_{\pi(1)},\ldots,a_nb_{\pi(n)}$ are (pairwise) distinct for
some permutation $\pi\in S_n$ (the symmetric group on
$\{1,\ldots,n\}$) if and only if $b_1\cdots b_n$ is the identity
element of $G$.

In 1999,  H. S. Snevily \cite{snevily} considered subsets with
cardinality $k$ of an abelian group $G$ (or simply $k$-subsets of
$G$), and proposed the following challenging conjecture.

\begin{Conjecture}\label{snev}{\rm (Snevily's Conjecture)}
Let $G$ be a multiplicatively written abelian group of odd order,
and let $A=\{a_1,\ldots, a_k\}$ and $B=\{b_1,\ldots, b_k\}$ be two
$k$-subsets of $G$. Then there is a permutation $\pi\in S_k$ such
that $a_1b_{\pi(1)},\ldots, a_kb_{\pi(k)}$ are distinct.
\end{Conjecture}

The above conjecture can be reformulated in terms of Latin
transversals of the Latin square formed by the Cayley
multiplication table of $G$. N. Alon \cite{alon} proved Snevily's
conjecture when $|G|$ is an odd prime, by using the Combinatorial
Nullstellensatz \cite{alon0}. In fact, Alon \cite{alon} obtained
the following stronger result.

\begin{Theorem}{\rm (Alon \cite{alon})}
Let $G$ be a cyclic group of prime order $p$. Let $k<p$ be a
positive integer. Let $A=\{a_1,\ldots ,a_k\}$ be a $k$-subset of
$G$ and $b_1,\ldots ,b_k$ be (not necessarily distinct) elements
of $G$. Then there is a permutation $\pi\in S_k$ such that
$a_1b_{\pi(1)},\ldots, a_kb_{\pi(k)}$ are distinct.
\end{Theorem}

In 2001, S. Dasgupta, G. K\'arolyi, O. Serra and B. Szegedy
\cite{das} proved Conjecture \ref{snev} for any cyclic group $G$
of odd order. Moreover, these authors extended Alon's result
(Theorem 1.2) to cyclic groups of prime power order as well as
elementary abelian groups.

\begin{Theorem}{\rm (Dasgupta-K\'arolyi-Serra-Szegedy \cite{das})}
Let $p$ be a prime and let $\alpha$ be a positive integer. Let $G$
be the cyclic group $C_{p^{\alpha}}$ or the elementary abelian
$p$-group $C_p^{\alpha}$. Assume that $A=\{a_1,\ldots ,a_k\}$ is a
$k$-subset of $G$ and $b_1,\ldots ,b_k$ are (not necessarily
distinct) elements of $G$, where $k$ is smaller than $p$. Then, for
some $\pi\in S_k$ the products $a_1b_{\pi(1)},\ldots, a_kb_{\pi(k)}$
are distinct.
\end{Theorem}

Motivated by Theorems 1.2 and 1.3, Dasgupta, et al. \cite[p.~23]{das}
proposed the following conjecture.

\begin{Conjecture}{\rm (The DKSS Conjecture)}\label{dasconj}
Let $G$ be a finite abelian group with $|G|>1$, and let $p(G)$ be
the smallest prime divisor of $|G|$. Let $k<p(G)$ be a positive
integer. Assume that $A=\{a_1,\ldots ,a_k\}$ is a $k$-subset of $G$
and $b_1,\ldots ,b_k$ are (not necessarily distinct) elements of
$G$. Then there is a permutation $\pi\in S_k$ such that
$a_1b_{\pi(1)},\ldots, a_kb_{\pi(k)}$ are distinct.
\end{Conjecture}

Let $G$ be a finite abelian group with $|G|>1$. When
$A=\{a_1,\ldots,a_k\}$ is a subgroup of $G$ with cardinality $k$ and
$b_1,\ldots ,b_k$ are (not necessarily distinct) elements of $A$, by
Hall's theorem, $a_1b_{\pi(1)},\ldots, a_kb_{\pi(k)}$ are distinct
for some $\pi\in S_k$ if and only if $b_1\cdots b_k=e$. As $G$
always has an element of order $p=p(G)$, which generates a cyclic
subgroup of order $p$, we see that the conclusion of Conjecture
\ref{dasconj} does not hold in general when $k=p(G)$. We also note
that Conjecture \ref{dasconj} implies Snevily's conjecture in the
case where $k<p(G)$. Using a group ring approach, W. D. Gao and D.
J. Wang \cite{gao} proved Conjecture~\ref{dasconj} for abelian
$p$-groups under the stronger assumption $k<\sqrt{2p}$. (See also [GGK].)

In this paper we confirm the DKSS conjecture under the extra
assumption that the second smallest prime divisor
of $|G|$ (if it exists) is greater than $k!$. It will be convenient to introduce
the following terminology.

\begin{Definition}\label{k-large} Let $k$ and $n>1$ be positive integers.
We say that $n$ is {\it $k$-large} if the smallest prime divisor of
$n$ is greater than $k$ and any other prime divisor of $n$ (if there
is any) is greater than $k!$. \end{Definition}

Here is our main result on the DKSS conjecture.

\begin{Theorem}\label{main} Let $G$ be a finite abelian group.
Let $A=\{a_1,\ldots, a_k\}$ be a $k$-subset of $G$, and $b_1,
\ldots,b_k$ be (not necessarily distinct) elements of $G$. Suppose
that either $A$ or $B=\{b_1,\ldots,b_k\}$ is contained in a subgroup
$H$ of $G$ and $|H|$ is $k$-large. Then there exists a permutation
$\pi\in S_k$ such that $a_1b_{\pi(1)},\ldots, a_kb_{\pi(k)}$ are
distinct.\end{Theorem}

Note that if $k$ is a positive integer and $p$ is a prime such that
$p>k$, then $p^{\alpha}$ is obviously $k$-large for every
$\alpha=1,2,3,\ldots$. Therefore we have the following immediate
corollary of Theorem~\ref{main}.

\begin{Corollary}\label{cor-0} Let $p$ be a prime. Assume that $G$ is an abelian $p$-group,
and $k$ is a positive integer such that $k<p$. Let $A=\{a_1,\ldots,
a_k\}$ be a $k$-subset of $G$, and $b_1,\ldots,b_k$ be (not
necessarily distinct) elements of $G$. Then there is a permutation
$\pi\in S_k$ such that $a_1b_{\pi(1)},\ldots, a_kb_{\pi(k)}$ are
distinct.\end{Corollary}

Obviously Theorem 1.6 implies that the DKSS conjecture is true for
$k=3$. By a case-by-case analysis we also can show that the DKSS
conjecture holds when $k=4$. Theorem~\ref{main} also implies that
the DKSS conjecture holds in the case where $A=\{a,a^2,\ldots,a^k\}$
with $a\in G\setminus\{e\}$.

\begin{Corollary}\label{cor-1} Let $G$ be a finite abelian group with
$|G|>1$. Let $a\not=e$ be an element of $G$ and let
$b_1,\ldots,b_k$ be (not necessarily distinct) elements of $G$.
Provided that $k<p(G)$, there is a permutation $\pi\in S_k$ such
that the products $a^ib_{\pi(i)}\ (i=1,\ldots,k)$ are
distinct.\end{Corollary}

It is interesting to compare Corollary \ref{cor-1} with the
following conjecture of Snevily \cite{snevily}: If $G=\langle
a\rangle$ is a cyclic group of order $n$, $k$ is a positive integer
less than $n$, and $b_1,\ldots,b_k$ are (not necessarily distinct)
elements of $G$, then there is a permutation $\pi\in S_k$ such that
$a^ib_{\pi(i)}$, $i=1,\ldots,k$, are distinct.

Our proof of Theorem~\ref{main} uses an exterior algebra approach.
For applications of symmetric product and alternating product
methods in combinatorics, we refer the reader to \cite[Chap.
6]{babaifrankl}. The exterior algebra approach can also be used to
give an alternative proof of the following theorem of Z. W. Sun
\cite{zwsun}.

\begin{Theorem}{\rm (Sun \cite{zwsun})}\label{minor}
Let $G$ be a finite cyclic group, and
let $A_1=\{a_{11},\ldots ,a_{1k}\}$, $A_2=\{a_{21},\ldots ,a_{2k}\}$,$\ldots$,
 $A_m=\{a_{m1},\ldots ,a_{mk}\}$ be $k$-subsets of $G$, where $m\in\{3,5,7,\ldots\}$.
 Then there exist permutations
 $\pi_i\in S_k$, $2\leq i\leq m$, such that $a_{1j}a_{2\pi_2(j)}\cdots a_{m\pi_m(j)}$, $j=1,\ldots ,k$, are all
 distinct.
\end{Theorem}

We remark that for $m=2,4,6,\ldots$, Sun could prove a similar
result under the assumption that all elements of $A_m$ have odd
order (cf. \cite{sun-JCTA} and \cite{zwsun}). This result in
particular implies that Snevily's conjecture is true for cyclic
groups of odd order, a result first proved by Dasgupta, et al. in
\cite{das}.

Recall that, for a matrix $M=(a_{ij})_{1\le i,j\le k}$ over a
commutative ring with identity, the determinant of $M$ and the
permanent of $M$ are defined by
$$\det(M)=\sum_{\sigma\in S_k}\sgn(\sigma)a_{1\sigma(1)}\cdots a_{k\sigma(k)}
\ \mbox{and}\ \per(M)=\sum_{\sigma\in S_k}a_{1\sigma(1)}\cdots
a_{k\sigma(k)}$$ respectively, where $\sgn(\sigma)$, the sign of
$\sigma\in S_k$, equals $1$ or $-1$ according as the permutation
$\sigma$ is even or odd.

To attack the Snevily conjecture via our approach, we propose the following conjecture.

\begin{Conjecture} \label{chi-det}
Let $G$ be a finite abelian group, and let $A=\{a_1,\ldots,a_k\}$ and $B=\{b_1,\ldots,b_k\}$
be two $k$-subsets of $G$. Let $K$ be any field
containing an element of multiplicative order $|G|$,
and let $\hat G$ be the character group of all group homomorphisms
from $G$ to $K^*=K\setminus\{0\}$. Then there are $\chi_1,\ldots,\chi_k\in\hat G$
such that  $\det(\chi_i(a_j))_{1\le i,j\le k}$
and $\det(\chi_i(b_j))_{1\le i,j\le k}$
are both nonzero.
\end{Conjecture}

When $G$ is cyclic, we may take $\chi_i=\chi^i$ for $i=1,\ldots ,k$,
where $\chi$ is a generator of ${\hat G}$. Then the two determinants
$\det(\chi_i(a_j))_{1\le i,j\le k}$ and $\det(\chi_i(b_j))_{1\le
i,j\le k}$ in the above conjecture are both nonzero since they are
Vandermonde determinants. Therefore we see that Conjecture
\ref{chi-det} is true for cyclic groups. We further mention that
when $G$ is a cyclic group of prime order and $K$ is the complex
field $\mathbb C$, for any distinct $\chi_1,\ldots,\chi_k\in\hat G$
and distinct $a_1,\ldots,a_k\in G$ we have $\det(\chi_i(a_j))_{1\le
i,j\le k}\not=0$ by the Chebotar\"ev theorem (cf. \cite{SL} and
\cite{tao}); this is stronger than what Conjecture~\ref{chi-det}
asserts. The general case of Conjecture~\ref{chi-det} seems to be
quite sophisticated.
\smallskip

We will show that Conjecture \ref{chi-det} holds when $A=B$ (see, Lemma \ref{lem 3.1}). Moreover, we have the following result.

\begin{Theorem} \label{chi-snev} {\rm (i)} Conjecture $\ref{chi-det}$ implies
Conjecture $\ref{snev}$.

{\rm (ii)} Let $G,A,B,\hat G$ be as in Conjecture $\ref{chi-det}$.
Assume that there is a $\pi\in S_k$ such that
$C=\{a_1b_{\pi(1)},\ldots,a_kb_{\pi(k)}\}$ is a $k$-set with
$\{a_1b_{\tau(1)},\ldots,a_kb_{\tau(k)}\}\not=C$ for all $\tau\in
S_k\setminus\{\pi\}$. Then there are $\chi_1,\ldots,\chi_k\in\hat G$
such that  $\det(\chi_i(a_j))_{1\le i,j\le k}$ and
$\det(\chi_i(b_j))_{1\le i,j\le k}$ are both nonzero. Also, there
are $\chi_1,\ldots,\chi_k\in\hat G$ such that
$\det(\chi_i(a_j))_{1\le i,j\le k}$ and $\per(\chi_i(b_j))_{1\le
i,j\le k}$ are both nonzero, and there are
$\psi_1,\ldots,\psi_k\in\hat G$ such that the permanents
$\per(\psi_i(a_j))_{1\le i,j\le k}$ and $\per(\psi_i(b_j))_{1\le
i,j\le k}$ are both nonzero.
\end{Theorem}

In the next section we will give some background on exterior
algebras and lay the basis for our new approach. (For readers who
are not familiar with exterior algebras, Northcott's book
\cite{north} is a good reference.) We are going to prove Theorem
\ref{main} and Corollary \ref{cor-1} in Section 3, and Theorems
\ref{minor} and \ref{chi-snev} in Section 4.

\section{An Auxiliary Proposition Motivated by Exterior Algebras}

 Let us recall some definitions and basic facts
related to exterior algebras.

Let $R$ be a commutative ring with identity, and let $M$ be a left
$R$-module. The {\it $n$th exterior power} of $M$, denoted by
$\bigwedge^nM$, comes equipped with an alternating multilinear map
$$M\times \cdots \times M\rightarrow {\bigwedge}^nM,\
\ (m_1,\ldots,m_n)\mapsto m_1\wedge \cdots \wedge m_n,$$ that is
universal: for an $R$-module $N$ and  an alternating multilinear map
$\beta:M\times \cdots \times M\rightarrow N$, there is a unique
linear map from $\bigwedge^nM$ to $N$ which takes $m_1\wedge \cdots
\wedge m_n$ to $\beta(m_1,\ldots ,m_n)$. Recall that a multilinear
map $\beta$ is {\it alternating} if $\beta(m_1,\ldots ,m_n)=0$
whenever two of the $m_i$ are equal. The exterior power
$\bigwedge^nM$ can be constructed as the quotient module of
$M^{\otimes n}$ (the $n$-th tensor power) by the submodule generated
by all those $m_1\otimes \cdots\otimes m_n$ with two of the $m_i$
equal. We naturally identify $\bigwedge^{0}M=R$ and
$\bigwedge^{1}M=M$. The {\it exterior algebra} of $M$, denoted by
$E(M)$, is the algebra $\oplus_{n\geq 0}\bigwedge^n(M)$, with
respect to the wedge product `$\wedge$'. This is a graded algebra.
By `graded' we mean that multiplying an element of $\bigwedge^mM$
with an element of $\bigwedge^nM$, gives an element of
$\bigwedge^{m+n}M$. A {\it skew derivation} on $E(M)$ is an
$R$-homomorphism $\Delta: E(M)\rightarrow E(M)$ such that
$$\Delta(xy)=(\Delta x)y+(-1)^n x(\Delta y),$$
for all $x\in \bigwedge^nM$ and $y\in E(M)$.

Next let $U$ be an $R$-module. Assume that we have a bilinear
mapping $\gamma: U\times M\rightarrow R$. Then, for any $u\in U$,
$\gamma(u,\cdot)$ is an $R$-module homomorphism from $M$ to $R$. By
\cite[Theorem 10, p.~96]{north}, there exists a unique skew
derivation $\Delta_u: E(M)\rightarrow E(M)$ that extends
$\gamma(u,\cdot)$. Furthermore, when $n>0$, $\Delta_u$ maps
$\bigwedge^nM$ into $\bigwedge^{n-1}M$, and it can be defined by
\[\Delta_u(m_1\wedge\cdots\wedge m_n)=
\sum_{i=1}^{n}(-1)^{i+1}\gamma(u,m_i)(m_1\wedge\cdots\wedge\widehat{m_i}\wedge\cdots\wedge
m_n),\] where $m_1,\ldots,m_n$ are arbitrary elements of $M$,
and $m_1\wedge\cdots\wedge\widehat{m_i}\wedge\cdots\wedge
m_n$ denotes the result of striking out $m_i$ from
$m_1\wedge\cdots\wedge m_i\wedge\cdots\wedge m_n$. For $u_1,\,u_2\in
U$, we can consider the composition $\Delta_{u_1}\circ\Delta_{u_2}$
of the $R$-module homomorphisms $\Delta_{u_1}$ and $\Delta_{u_2}$ in
the usual sense. We are now ready to state the following result from
\cite{north}.

\begin{Lemma}{\rm (\cite[Corollary, p.~100]{north})}\label{lem 1} Using the above
notation, for $u_1,\ldots,u_k\in U$ and $m_1,\ldots,m_k\in M$, we
have
$$(\Delta_{u_1}\circ\cdots\circ\Delta_{u_k})(m_1\wedge\cdots\wedge
m_k)=(-1)^{k(k-1)/2}\det\left(\gamma(u_i,m_j)\right)_{1\le i,j\le k}.$$
\end{Lemma}

Lemma \ref{lem 1} leads us to the following useful proposition.

\begin{Proposition}\label{prop} Let $G$ be a finite abelian group. Let $\hat{G}$ denote the
group of characters from $G$ to $K^*=K\setminus\{0\}$, where $K$ is
a field. Let $a_1,\ldots,a_k,b_1,\ldots,b_k\in G$ and
$\chi_1,\ldots,\chi_k\in\hat G$. Suppose that both $\det(M_a)$ and
$\per(M_b)$ are nonzero, where $M_a=(\chi_i(a_j))_{1\le i,j\le k}$
and $M_b=(\chi_i(b_j))_{1\le i,j\le k}$. Then there is a permutation
$\pi\in S_k$ such that the products
$a_1b_{\pi(1)},\ldots,a_kb_{\pi(k)}$ are distinct.
\end{Proposition}
\noindent{\it Proof}. For the purpose of applying Lemma~\ref{lem 1},
we set $R:=K$, $M:=KG$ (the group algebra of $G$ over $K$), and
$U:=K\hat{G}$ (the group algebra of $\hat{G}$ over $K$). The mapping
$\gamma: U\times M\rightarrow K$ is defined as follows: First define
$\gamma: \hat{G}\times G\rightarrow K$ by setting $\gamma(\chi,
g):=\chi(g)$ for $\chi\in \hat{G}$ and $g\in G$; next we extend
$\gamma$ bilinearly to a map from $U\times M$ to $K$. The resulting
map is still denoted by $\gamma$ and it is bilinear.

 For any $\pi\in S_k$ we set
$$Q_{\pi}:=a_1b_{\pi(1)}\wedge\cdots\wedge a_kb_{\pi(k)}\in{\bigwedge}^kM.$$
Let $M_{a,b}^{\pi}$ be the $k\times k$ matrix with $(i,j)$-entry
equal to $\chi_i(a_jb_{\pi(j)})$. By Lemma \ref{lem 1}, we have
\begin{equation}\label{preequ}\aligned
(\Delta_{\chi_1}\circ\cdots\circ\Delta_{\chi_k})(Q_{\pi})
=&(-1)^{k(k-1)/2}\det(M_{a,b}^{\pi})\\
=&(-1)^{k(k-1)/2}\sum_{\sigma\in
S_k}\sgn(\sigma)\chi_1(a_{\sigma(1)}b_{\pi\sigma(1)})\cdots\chi_k(a_{\sigma(k)}b_{\pi\sigma(k)}).
\endaligned
\end{equation}
Summing (\ref{preequ}) over $\pi\in S_k$, we obtain
\begin{equation}\label{mainequ}\aligned
&(\Delta_{\chi_1}\circ\cdots\circ\Delta_{\chi_k})\bigg(\sum_{\pi\in S_k}Q_{\pi}\bigg)
 \\=&(-1)^{k(k-1)/2}\sum_{\sigma\in S_k}\sum_{\pi\in S_k}
 \chi_1(a_{\sigma(1)}b_{\pi\sigma(1)})\cdots\chi_k(a_{\sigma(k)}b_{\pi\sigma(k)})\sgn(\sigma)
 \\=&(-1)^{k(k-1)/2}\sum_{\sigma\in S_k}\sgn(\sigma)\chi_1(a_{\sigma(1)})\cdots\chi_k(a_{\sigma(k)})
 \sum_{\pi\in S_k} \chi_1(b_{\pi(1)})\cdots\chi_k(b_{\pi(k)})\\
=&(-1)^{k(k-1)/2}\det(M_a){\rm per}(M_b).
\endaligned
\end{equation}
As $\det(M_a){\rm per}(M_b)\not=0$, we have
 $\sum_{\pi\in S_k}Q_{\pi}\neq 0$. So there exists a $\pi\in S_k$ such that
$a_1b_{\pi(1)}\wedge\cdots\wedge a_kb_{\pi(k)}$ is nonzero,
implying that $a_1b_{\pi(1)},\ldots, a_kb_{\pi(k)}$ are
distinct. The proof is now complete.\\

\section{Proofs of Theorem \ref{main} and Corollary \ref{cor-1}}

\begin{Lemma} \label{lem 3.1} Let $G$ be a finite abelian group
and let $K$ be a field containing an element of multiplicative order
$|G|$. Let $\hat G$ be the group of characters from $G$
to $K^*$, and let $a_1,\ldots,a_k\in G$. Then $a_1,\ldots,a_k$ are
distinct if and only if there are (distinct)
$\chi_1,\ldots,\chi_k\in\hat G$ such that $\det(\chi_i(a_j))_{1\le
i,j\le k}\neq0$.
Also, there exist $\chi_1,\ldots,\chi_k\in\hat G$ with
$\per(\chi_i(a_j))_{1\le i,j\le k}\not=0$ provided that $a_1,\ldots,a_k$ are distinct.
\end{Lemma}
\noindent{\it Proof}. If $a_s=a_t$ for
some $1\le s<t\le k$, then for any $\chi_1,\ldots,\chi_k\in\hat G$
the determinant $\det(\chi_i(a_j))_{1\le i,j\le k}$ vanishes since
the $s$th column and $t$th column of the matrix $(\chi_i(a_j))_{1\le
i,j\le k}$ are identical.

Now suppose that $a_1,\ldots,a_k$ are distinct.
If the characteristic of $K$ is a prime $p$ dividing $|G|$, then
$$(x^{|G|/p}-1)^p=x^{|G|}-1\quad\mbox{for all}\ x\in K,$$
which contradicts the assumption that $K$ contains an element of
multiplicative order $|G|$. So we have $|G|1\not=0$, where $1$ is
the identity of the field $K$. It is well known that
$$\sum_{\chi\in \hat G}\chi(a)=\begin{cases}0, &\mbox{if}\ a\in G\setminus\{e\},
\\|G|1, &\mbox{if}\ a=e.\end{cases}$$

Observe that
\begin{align*} &\sum_{\chi_1,\ldots,\chi_k\in\hat
G}\chi_1(a_1^{-1})\cdots\chi_k(a_k^{-1})\det(\chi_i(a_j))_{1\le
i,j\le k} \\=&\sum_{\chi_1,\ldots,\chi_k\in\hat
G}\chi_1(a_1^{-1})\cdots\chi_k(a_k^{-1})\sum_{\pi\in
S_k}\sgn(\pi)\prod_{i=1}^k\chi_i(a_{\pi(i)})
\\=&\sum_{\chi_1,\ldots,\chi_k\in\hat G}\ \sum_{\pi\in S_k}\sgn(\pi)\prod_{i=1}^k\chi_i(a_{\pi(i)}a_i^{-1})
\\=&\sum_{\pi\in S_k}\sgn(\pi)\prod_{i=1}^k\sum_{\chi_i\in\hat G}\chi_i(a_{\pi(i)}a_i^{-1})
\\=&\sgn(I)\prod_{i=1}^k(|G|1)=(|G|1)^k\not=0,
\end{align*}
where $I$ is the identity permutation in $S_k$.
So there are $\chi_1,\ldots,\chi_k\in\hat G$ such that
$\det(\chi_i(a_j))_{1\le i,j\le k} \not=0$ and hence $\chi_1,\ldots,\chi_k$
are distinct. Similarly, by removing $\sgn(\pi)$ in the above calculations, we see that there exist $\chi_1,\ldots,\chi_k\in\hat G$ such that
$\per(\chi_i(a_j))_{1\le i,j\le k} \not=0$.
This concludes the proof.

\begin{Remark}\label{rem-1} If we apply Lemma \ref{lem 3.1} with $k=|G|$ then we obtain the following classical result:
The matrix $T=\big(\chi(g)\big)_{\chi\in \hat{G},\,g\in G }$ is
nonsingular; in other words, all the characters in $\hat G$ are
linearly independent over the field $K$. It is well known that  all
the characters in $\hat G$ actually form a basis of the vector space
$$K^G=\{f:\ f\ \mbox{is a function from}\ G\ \mbox{to}\ K\}$$
over the field $K$.
\end{Remark}

\begin{Lemma} {\rm (Sun \cite[Lemma 3.1]{sun03})} \label{lem-2}  Let $\lambda_1,\ldots,\lambda_k$ be
complex $n$th roots of unity. Suppose that
$c_1\lambda_1+\cdots+c_k\lambda_k=0$, where $c_1,\ldots,c_k$ are
nonnegative integers. Then $c_1+\cdots+c_k$ can be written in the
form $\sum_{p|n}px_p$, where the sum is over all prime divisors of
$n$ and the $x_p$ are nonnegative integers.
\end{Lemma}

 This lemma appeared in \cite{sun03} explicitly, and it follows from Lemma 9 in Sun \cite{sun96}.

\medskip

We are now in a position to prove
Theorem~\ref{main}.\\

\noindent\textsl{Proof of Theorem \ref{main}.} Let $\hat G$ be the
group of all complex-valued characters of $G$, and let
$p_1<\cdots<p_r$ be all the distinct prime divisors of $|H|$. If
$k!=p_1x_1+\cdots+p_rx_r$ for some nonnegative integers $x_1,\ldots,x_r$,
then $k!=p_ix_i$ for some $1\le i\le r$ (since $p_{i_1}+p_{i_2}>k!$
for all $1\le i_1<i_2\le r$), thus $p_i$ divides $k!$, contradicting
the fact $k<p(H)\le p_i$. In view of this and Lemma \ref{lem-2}, we
see that if $\zeta_{\pi}\in\mathbb C$ and $\zeta_{\pi}^{|H|}=1$ for all $\pi\in S_k$ then
$\sum_{\pi\in S_k}\zeta_{\pi}\not=0$.

\medskip

{\it Case} 1. $B=\{b_1,\ldots,b_k\}\subseteq H$.

In this case $b_1^{|H|}=\cdots=b_k^{|H|}=e$. As $a_1,\ldots,a_k$
are distinct, by Lemma \ref{lem 3.1} there are
$\chi_1,\ldots,\chi_k\in\hat G$ such that $\det(\chi_i(a_j))_{1\le
i,j\le k}\not=0$. Observe that
$$\per(\chi_i(b_j))_{1\le i,j\le k}=\sum_{\pi\in S_k}\prod_{i=1}^k\chi_i(b_{\pi(i)})$$
and
$$\bigg(\prod_{i=1}^k\chi_i(b_{\pi(i)})\bigg)^{|H|}=\prod_{i=1}^k\chi_i(b_{\pi(i)}^{|H|})=\prod_{i=1}^k\chi_i(e)=1.$$
By the previous discussion,  $\per(\chi_i(b_j))_{1\le i,j\le
k}\not=0$. It follows from Proposition \ref{prop} that
$a_1b_{\pi(1)},\ldots,a_kb_{\pi(k)}$ are distinct for some $\pi\in
S_k$.

\medskip

{\it Case} 2. $A\subseteq H$.

 Suppose that $p_1,\ldots,p_r,\ldots,p_s\ (s\ge r)$ are all the distinct prime divisors of $|G|$.
 It is well known that $G=P_1\cdots P_r\cdots P_s$, where each $P_i=\mbox{Syl}_{p_i}(G)$ is the unique Sylow $p_i$-subgroup of $G$.
 For each $i=1,\ldots,r$, we have $\mbox{Syl}_{p_i}(H)\subseteq P_i$ by Sylow's second theorem. So
 $$H=\mbox{Syl}_{p_1}(H)\cdots\mbox{Syl}_{p_r}(H)\subseteq H_1:=P_1\cdots P_r.$$
As $a_1,\ldots,a_k$ are distinct elements of $H_1$, by Lemma
\ref{lem 3.1} there are $\chi_1,\ldots,\chi_k\in\widehat {H_1}$ such that
$\det(\chi_i(a_j))_{1\le i,j\le k}\not=0$. For each $i=1,\ldots,k$
we can extend $\chi_i$ to a character of $G$ by setting
$\chi_i(h_1h_2)=\chi_i(h_1)$ for all $h_1\in H_1$ and $h_2\in H_2$,
where $H_2=P_{r+1}\cdots P_s$ if $r<s$, and $H_2=\{e\}$ if $r=s$. It
follows that $\chi_i(g)^{|H_1|}=1$ for all $g\in G=H_1H_2$. Thus
$$\bigg(\prod_{i=1}^k\chi_i(b_{\pi(i)})\bigg)^{|H_1|}=1\quad\mbox{for all}\;\pi\in S_k.$$
Note that $|H_1|$ is $k$-large since $|H|$ is. As in Case 1, we get
$\per(\chi_i(b_j))_{1\le i,j\le k}\not=0$ and hence the
desired result follows.
\medskip

Combining the above, we have completed the proof of Theorem \ref{main}.\\

\noindent{\it Proof of Corollary \ref{cor-1}}. Let $p$ be any prime
divisor of the order of $a$, and let $P={\rm Syl}_p(G)$, the unique
Sylow $p$-subgroup of $G$. Since the order of $a$ is a multiple of
$p$ and $|G|/|P$ is relatively prime to $p$, the order of $a$ does
not divide $|G|/|P|$. Therefore $\tilde a=a^{|G|/|P|}\not=e$. As
$(\tilde a)^{|P|}=a^{|G|}=e$, the group generated by $\tilde a$ is a
$p$-subgroup of $G$. So we have $\tilde a \in P\setminus\{e\}$.

The condition $k<p(G)$ ensures that $1,\ldots,k$ are pairwise
incongruent modulo $p$. Our following argument actually yields a
refinement of the stated result.

Let $i_1,\ldots,i_k$ be any integers pairwise incongruent modulo
$p$. Then $\tilde a^{i_1},\ldots,\tilde a^{i_k}$ are distinct
elements of $P$. Set ${\tilde b_i}=b_i^{|G|/|P|}$. By Theorem
\ref{main}, there is a permutation $\pi\in S_k$ such that all the
elements
$$\tilde a^{i_j}\widetilde {b_{\pi(j)}}=(a^{i_j}b_{\pi(j)})^{|G|/|P|}\ \ (j=1,\ldots,k)$$
are distinct. It follows that
$a^{i_1}b_{\pi(1)},\ldots,a^{i_k}b_{\pi(k)}$ are distinct. We are
done.

\section{Proofs of Theorems \ref{minor} and \ref{chi-snev}}

First, we prove Theorem~\ref{minor} via the exterior algebra approach.\\

\noindent\textsl{Proof of Theorem \ref{minor}.} We choose $K=\mathbb{C}$ (the field
of complex numbers), and work with the group algebra $M:=KG$. For
$\pi_2,\ldots,\pi_m\in S_k$ we define
$$Q_{\pi_2,\ldots,\pi_m}=\prod_{i=2}^m{\rm sgn}(\pi_i)
\left(a_{11}a_{2\pi_2(1)}\cdots a_{m\pi_m(1)}\wedge\cdots \wedge
     a_{1k}a_{2\pi_2(k)}\cdots a_{m\pi_m(k)}\right)\in{\bigwedge}^kM.$$
Let $\chi_1,\ldots ,\chi_k$ be distinct (complex) characters of $G$,
and set $M_{i}=(\chi_{\ell}(a_{ij}))_{1\leq\ell,j\leq k}$ for
$i=1,2,\ldots ,m$. Then, by similar computations to those in
(\ref{preequ}) and (\ref{mainequ}) (paying attention to the signs
involved and noting that $m$ is odd), we get
$$\Delta_{\chi_1}\circ\cdots\circ\Delta_{\chi_k}
\bigg(\sum_{\pi_2,\ldots ,\pi_m\in S_k}Q_{\pi_2,\ldots
,\pi_m}\bigg)=(-1)^{k(k-1)/2}\prod_{i=1}^m\det(M_{i}).$$ Since $G$
is cyclic, we have ${\hat G}=\{1,\chi,\chi^2,\ldots\}$, where $\chi$
is a generator of ${\hat G}$. If $\chi(a_{is})=\chi(a_{it})$, then
$\sum_{\psi\in\hat G}\psi(a_{is}a_{it}^{-1})=|G|1\not=0$, hence
$a_{is}a_{it}^{-1}=e$ and thus $s=t$ (since $a_{i1},\ldots,a_{ik}$
are distinct). Therefore, if we choose $\chi_{\ell}=\chi^{\ell-1}$
for $\ell=1,\ldots ,k$, then $\det(M_{i})\neq 0$ for all $1\leq
i\leq m$ since each $M_{i}$ is a Vandermonde matrix whose second row
is $(\chi(a_{ij}))_{1\leq j\leq k}$. With the above choice of
$\chi_{\ell}$, $1\leq \ell\leq k$, we have
$$\Delta_{\chi_1}\circ\cdots\circ\Delta_{\chi_k}
\bigg(\sum_{\pi_2,\ldots ,\pi_m\in S_k}Q_{\pi_2,\ldots
,\pi_m}\bigg)\neq 0.$$ Hence there exist $\pi_2,\ldots ,\pi_m\in
S_k$ such that $Q_{\pi_2,\ldots,\pi_m}\neq 0$, implying that
$a_{1j}a_{2\pi_2(j)}\cdots a_{m\pi_m(j)}$, $1\leq j\leq k$, are all
distinct. \\

Next, we give a proof of Theorem \ref{chi-snev}.
\medskip

\noindent\textsl{Proof of Theorem \ref{chi-snev}.}
(i) We want to prove Conjecture \ref{snev} under the assumption that Conjecture \ref{chi-det} holds.

As $|G|$ is odd, we have $2^{\varphi(|G|)}\equiv1\ (\mbox{mod}\
|G|)$ by Euler's theorem, where $\varphi$ is the Euler totient
function. Let $K$ be the finite field
$\mathbb{F}_{2^{\varphi(|G|)}}$. Then the cyclic group
$K^*=K\setminus\{0\}$ has an element of multiplicative order
$|G|$. Let ${\hat G}$ be the group of characters from $G$ to
$K^*$. By Conjecture \ref{chi-det} there are
$\chi_1,\ldots,\chi_k\in\hat G$ such that $\det(M_a)$ and $\det(M_b)$ are both nonzero,
where $M_a=(\chi_i(a_j))_{1\le i,j\le k}$
and $M_b=(\chi_i(b_j))_{1\le i,j\le k}$. As $K$ is of
characteristic 2, we have $\per (M_b)=\det (M_b)\not=0$. Applying
Proposition \ref{prop} we obtain that
$a_1b_{\pi(1)},\ldots,a_kb_{\pi(k)}$ are distinct for some $\pi\in
S_k$.

(ii) Set $$\Sigma:=\sum_{\chi_1,\ldots,\chi_k\in\hat G}
\chi_1^{-1}(a_1b_{\pi(1)})\cdots\chi_k^{-1}(a_kb_{\pi(k)})\det((\chi_i(a_j))_{1\le
i,j\le k} \det((\chi_i(b_j))_{1\le i,j\le k}.$$

Observe that
\begin{align*}\Sigma=&\sum_{\chi_1,\ldots,\chi_k\in\hat G}\sum_{\sigma\in S_k}\sgn(\sigma)
\prod_{i=1}^k\chi_i(a_{\sigma(i)}a_i^{-1})
\sum_{\tau\in S_k}\sgn(\tau\sigma)\prod_{i=1}^k\chi_i(b_{\tau\sigma(i)}b_{\pi(i)}^{-1})
\\=&\sum_{\tau\in S_k}\sgn(\tau)\sum_{\sigma\in S_k}\prod_{i=1}^k
\sum_{\chi_i\in\hat G}\chi_i(a_{\sigma(i)}b_{\tau\sigma(i)}(a_ib_{\pi(i)})^{-1})
\\=&\sum_{\tau\in S_k}\sgn(\tau)\sum_{\sigma\in S_k\ \atop{a_{\sigma(i)}b_{\tau\sigma(i)}=a_ib_{\pi(i)}
\ {\rm for}\ i=1,\ldots,k}}(|G|1)^k.
\end{align*}
If $\tau\in S_k$ and $\{a_1b_{\tau(1)},\ldots,a_kb_{\tau(k)}\}=C$, then
 $\tau$ is identical to $\pi$ by the assumption of Theorem~\ref{chi-snev}(ii), and the identity permutation $I$
 is the unique $\sigma\in S_k$ such that $a_{\sigma(i)}b_{\tau\sigma(i)}=a_ib_{\pi(i)}$ for all
$i=1,\ldots,k$. Therefore
$$\Sigma=\sgn(\pi)(|G|1)^k\not=0.$$
So there are $\chi_1,\ldots,\chi_k\in \hat G$ such that
$$\det((\chi_i(a_j))_{1\le i,j\le k}\det((\chi_i(b_j))_{1\le i,j\le k}\not=0$$
and hence neither $\det((\chi_i(a_j))_{1\le i,j\le k}$ nor
$\det((\chi_i(b_{j}))_{1\le i,j\le k}$ vanishes.

The remaining results in Theorem \ref{chi-snev}(ii) can be easily proved in a similar way.

\end{document}